\documentclass[10pt]{amsart}
\usepackage[T1]{fontenc}
\usepackage{lmodern}
\usepackage{geometry}
\usepackage[latin1] {inputenc}
\usepackage{amsmath}
\usepackage{amsfonts, amssymb, textcomp}
\usepackage[colorlinks=flase, linkcolor=red,urlcolor=green, citecolor=blue]{hyperref}
\usepackage{subeqnarray}

\usepackage{latexsym}
\usepackage{fancyhdr}
\usepackage{longtable}
\usepackage{amsmath, amssymb}
\usepackage{graphicx}

\numberwithin{equation}{section}

\theoremstyle{plain}
\newtheorem{theorem}{Theorem}[section]
\newtheorem{lemma}[theorem]{Lemma}

\theoremstyle{definition}
\newtheorem{definition}[theorem]{Definition}
\newtheorem{example}[theorem]{Example}
\newtheorem{remark}[theorem]{Remark}

\newcommand{\RR}{\mathbb{R}}

\newcommand{\NN}{\mathbb{N}}

\let\on=\operatorname

\newenvironment{demo}[1]{\par\smallskip\noindent{\bf #1.}}{\par\smallskip}

\makeatletter
\@namedef{subjclassname@2020}{%
  \textup{2020} Mathematics Subject Classification}
\makeatother

\title[On the class of almost subadditive weight functions]
{On the class of almost subadditive weight functions}

\author[G.~Schindl]{Gerhard Schindl}

\address{G.~Schindl: Fakult\"at f\"ur Mathematik, Universit\"at Wien, Oskar-Morgenstern-Platz~1, A-1090 Wien, Austria.}
\email{gerhard.schindl@univie.ac.at}

\begin{document}

\begin{abstract}
We answer a question from A. V. Abanin and P. T. Tien about so-called almost subadditive weight functions in the sense of Braun-Meise-Taylor. Using recent knowledge of a growth index for functions, crucially appearing in the ultraholomorphic setting, we are able to show the existence of weights such that there does not exist an equivalent almost subadditive function and hence almost subadditive weights are a proper subclass.
\end{abstract}

\thanks{This research was funded in whole by the Austrian Science Fund (FWF) project 10.55776/P33417}
\date{\today}

\maketitle

\section{Introduction}
In \cite{Abanin10} it has been shown that {\itshape almost subadditive weight functions} are sufficient to describe the theory of ultradistributions in the sense of Braun-Meise-Taylor; see \cite{BraunMeiseTaylor90}. Thus the authors are dealing with non-quasianalytic weight functions $\omega$; see Section \ref{growthsection} for precise definitions and conditions. All relevant information about growth and regularity properties for weight functions and weight sequences can be found in Section \ref{weightgrowthsect}.

In their closing discussion in \cite[Rem. 1.5]{Abanin10} the authors have asked the following:\vspace{6pt}

Let $\omega$ be a given non-quasianalytic weight function (in the sense of Braun-Meise-Taylor), does there always exist a weight $\sigma$ such that
\begin{itemize}
\item[$(*)$] $\sigma$ is {\itshape equivalent} to $\omega$ which means that $0<\liminf_{t\rightarrow+\infty}\frac{\omega(t)}{\sigma(t)}\le\limsup_{t\rightarrow+\infty}\frac{\omega(t)}{\sigma(t)}<+\infty$, and such that
\item[$(*)$] $\sigma$ is {\itshape almost subadditive:}
$$\forall\;q>1\;\exists\;C_q\ge 1\;\forall\;s,t\ge 0:\;\;\;\sigma(s+t)\le q(\sigma(s)+\sigma(t))+C_q.$$
\end{itemize}
The notion of {\itshape almost subadditivity} is more recent and more general than the classical (standard) assumption {\itshape subadditivity} appearing in the Beurling-Bj\"{o}rck framework; see e.g. \cite{Bjorck66}. On the other hand we mention that in the literature different ''subadditivity-like'' conditions for (associated) weight functions have been studied and used; see e.g. \cite{subaddlike} and the citations there.\vspace{6pt}

The aim of this article is to answer this question and to prove that in general there does not exist an equivalent almost subadditive weight. This fact shows that the class of almost subadditive weight functions is a proper subclass and hence the main result \cite[Thm. 1.4]{Abanin10} means that this particular subclass of weight functions is sufficient to describe the theory of ultradistributions in the sense of Braun-Meise-Taylor. Our statements can be found in Section \ref{Abaninmainressect}. In this work, since we are mainly focusing on functions and their growth properties, it is not required to introduce the associated ultradifferentiable function spaces explicitly; we refer to \cite{BraunMeiseTaylor90}, \cite{BonetMeiseMelikhov07}, \cite{Komatsu73}, \cite{compositionpaper}.\vspace{6pt}

In order to treat this problem, crucially we involve recent knowledge about the growth index $\gamma(\omega)$ from \cite{index}; see Section \ref{growthindexsect} for precise definitions and explanations. This value becomes relevant when proving extension results in the ultraholomorphic setting; see \cite{sectorialextensions}, \cite{sectorialextensions1}. We establish a connection between the notion of almost subadditivity and $\gamma(\omega)$ (see Lemma \ref{crucialindexlemma} and Remark \ref{crucialindexlemmarem}). Involving this information, intrinsic properties of $\gamma(\omega)$, in particular using the fact that this growth index is preserved under equivalence of weights, we show the main result Theorem \ref{Abaninmainthm}.

Note that using this method we are able to give more information: On the one hand, we are working within a more general setting of weight functions; see Definition \ref{weightfctdef}. Thus our basic assumptions on $\omega$ are weaker than being a weight in the sense of Braun-Meise-Taylor, but also weaker than the standard requirements assumed in \cite{Bjorck66} and even slightly weaker than in \cite{sectorialextensions}, \cite{sectorialextensions1}. This is possible by the general setting studied in \cite[Sect. 2.3 \& 2.4]{index}. On the other hand, instead of almost subadditivity, we treat a  more general condition (see \eqref{almostsubaddgeneral}) and obtain that the above problem is not depending on the notion of (non-)quasianalyticity. So the main result is valid for quasianalytic weights as well; see Remark \ref{Abaninmainthmrem} for more detailed comments.\vspace{6pt}

However, Theorem \ref{Abaninmainthm} provides abstract information and the aim in Theorem \ref{Abaninexamplethm} and Theorem \ref{Abaninexamplethm1} is to construct explicitly weights with the desired property. There we are focusing on so-called {\itshape associated weight functions} $\omega_M$ with $M\in\RR_{>0}^{\NN}$ being a {\itshape weight sequence;} see Section \ref{assofuncsection}. The advantage in this case is that the desired growth properties for $\omega_M$ can be expressed in terms of $M$ and so one is able to proceed by constructing an appropriate weight sequence $M$. Such explicit constructions of (counter-)examples can be useful for further applications in different contexts. Here, the growth index $\gamma(M)$ originally introduced by V. Thilliez in \cite{Thilliezdivision} is becoming relevant; see again \cite{index} for more details and its relation to $\gamma(\omega_M)$.



\section{Weights, growth conditions and the index $\gamma(\omega)$}\label{weightgrowthsect}

\subsection{General notation}
We use the notation $\RR_{>0}:=(0,+\infty)$, $\NN:=\{0,1,2,\dots\}$ and $\NN_{>0}:=\{1,2,\dots\}$.

\begin{definition}\label{weightfctdef}
A function $\omega:[0,+\infty)\rightarrow[0,+\infty)$ is called a {\itshape weight function} if
\begin{itemize}
\item[$(*)$] $\omega$ is non-decreasing and

\item[$(*)$] $\lim_{t\rightarrow+\infty}\omega(t)=+\infty$.
\end{itemize}
\end{definition}
These are the basic assumptions resp. the notion for being a weight which is used in \cite[Sect. 2.3 \& 2.4]{index}. These requirements are sufficient to introduce and work with the crucial index $\gamma(\omega)$; see Section \ref{growthindexsect} for details.

\subsection{Relevant growth conditions}\label{growthsection}
We list now conditions for weight functions $\omega:[0,+\infty)\rightarrow[0,+\infty)$; they are known and standard in the ultradifferentiable setting involving weight functions in the sense of {\itshape Braun-Meise-Taylor;} see \cite{BraunMeiseTaylor90}. We use the abbreviations appearing in \cite{dissertation}.

\begin{itemize}
\item[\hypertarget{om1}{$(\omega_1)}$] $\omega(2t)=O(\omega(t))$ as $t\rightarrow+\infty$.
	
\item[\hypertarget{om2}{$(\omega_2)$}] $\omega(t)=O(t)$ as $t\rightarrow+\infty$.
	
\item[\hypertarget{om3}{$(\omega_3)$}] $\log(t)=o(\omega(t))$ as $t\rightarrow+\infty$.
	
\item[\hypertarget{om4}{$(\omega_4)$}] $\varphi_{\omega}:t\mapsto\omega(e^t)$ is a convex function (on $\RR$).
	
\item[\hypertarget{om5}{$(\omega_5)$}] $\omega(t)=o(t)$ as $t\rightarrow+\infty$.
	
\item[\hypertarget{omnq}{$(\omega_{\text{nq}})$}] $\int_1^{\infty}\frac{\omega(t)}{t^2}dt<+\infty.$
	
\item[\hypertarget{omsnq}{$(\omega_{\text{snq}})$}] $\exists\;C>0\;\forall\;y>0: \int_1^{\infty}\frac{\omega(y t)}{t^2}dt\le C\omega(y)+C$.
\end{itemize}

In the literature, when being defined on $\RR^d$, usually $\omega$ is radially extended, i.e. $\omega(x):=\omega(|x|)$ for any $x\in\RR^d$.

We call a weight function $\omega$ {\itshape normalized} if $\omega(t)=0$ for all $0\le t\le 1$, {\itshape non-quasianalytic} if \hyperlink{omnq}{$(\omega_{\on{nq}})$} is valid and {\itshape quasianalytic} if \hyperlink{omnq}{$(\omega_{\on{nq}})$} is violated.\vspace{6pt}

\begin{definition}
A function $\omega:[0,+\infty)\rightarrow[0,+\infty)$ is called a BMT-weight function (BMT standing for Braun-Meise-Taylor), if $\omega$ is continuous, non-decreasing, $\omega(0)=0$, and such that $\omega$ satisfies \hyperlink{om1}{$(\omega_1)$}, \hyperlink{om3}{$(\omega_3)$} and \hyperlink{om4}{$(\omega_4)$}.

If $\omega$ is in addition normalized, then it is called a normalized BMT-weight function; if in addition \hyperlink{omnq}{$(\omega_{\on{nq}})$} is valid, then $\omega$ is called a non-quasianalytic BMT-weight function and it is called quasianalytic, if \hyperlink{omnq}{$(\omega_{\on{nq}})$} fails.
\end{definition}

In \cite{Abanin10} it has been assumed that $\omega$ is continuous and non-decreasing. Moreover, the following abbreviations have been used: $(\alpha)$ is \hyperlink{om1}{$(\omega_1)$}, $(\beta)$ is \hyperlink{omnq}{$(\omega_{\on{nq}})$}, $(\gamma)$ denotes \hyperlink{om3}{$(\omega_3)$} and $(\delta)$ is \hyperlink{om4}{$(\omega_4)$}. Summarizing, in \cite{Abanin10} non-quasianalytic BMT-weight functions (not necessarily satisfying $\omega(0)=0$, but this failure is not effecting the results) have been considered and this set is denoted by $\Omega$ there. Note that in the proof of the main technical statement \cite[Prop. 1.3]{Abanin10} non-quasianalyticity has been applied crucially.\vspace{6pt}

Let $\sigma,\tau: [0,+\infty)\rightarrow[0,+\infty)$ be arbitrary, then we write $\sigma\hypertarget{ompreceq}{\preceq}\tau$ if
\begin{equation*}\label{bigOrelation}
\tau(t)=O(\sigma(t))\;\text{as}\;t\rightarrow+\infty,	
\end{equation*}
and $\sigma$ and $\tau$ are called {\itshape equivalent,} written $\sigma\hypertarget{sim}{\sim}\tau$, if
$$\sigma\hyperlink{ompreceq}{\preceq}\tau\;\text{and}\;\tau\hyperlink{ompreceq}{\preceq}\sigma.$$
This relation is crucial in the ultradifferentiable weight function setting (in the sense of Braun-Meise-Taylor), see \cite[Lemma 5.16, Cor. 5.17]{compositionpaper}, and it is precisely the same notion of equivalence used in \cite{Abanin10}; see \cite[Sect. 1, p. 297 \& Rem. 1.5]{Abanin10}. Moreover, it also has been considered in \cite{index}; see \cite[Sect. 2.1 \& Rem. 2.6 $(5)$]{index}. It is immediate that all properties listed above, except the convexity condition \hyperlink{om4}{$(\omega_4)$}, are preserved under equivalence and that normalization can be assumed w.l.o.g. by switching to an equivalent weight.\vspace{6pt}

According to \cite[Def. 1.1]{Abanin10} we say that a weight function $\omega:[0,+\infty)\rightarrow[0,+\infty)$ is {\itshape almost subadditive} if
\begin{equation}\label{almostsubadd}
\forall\;q>1\;\exists\;C_q\ge 1\;\forall\;s,t\ge 0:\;\;\;\omega(s+t)\le q(\omega(s)+\omega(t))+C_q.
\end{equation}

This should be compared with the more original condition {\itshape subadditivity} in the Beurling-Bj\"{o}rck framework (see \cite{Bjorck66}):

\begin{itemize}
\item[\hypertarget{sub}{$(\omega_{\text{sub}})$}] $\omega(s+t)\le\omega(s)+\omega(t)$ for all $s,t\ge 0$.
\end{itemize}

Obviously \hyperlink{sub}{$(\omega_{\on{sub}})$} implies \eqref{almostsubadd} with the uniform choice $C_q=0$ for all $q>1$. {\itshape Note:} If the inequality in \eqref{almostsubadd} is valid for some $q_0>1$, then for all $q\ge q_0$ as well when choosing $C_q=C_{q_0}$. So, in order to verify almost subadditivity, one has to study $q\rightarrow 1$. Note that by definition neither \hyperlink{sub}{$(\omega_{\on{sub}})$} nor \eqref{almostsubadd} is preserved under equivalence of weight functions in general.

\subsection{The growth index $\gamma(\omega)$}\label{growthindexsect}
We briefly recall the definition of the growth index $\gamma(\omega)$ and some consequences; see \cite[Sect. 2.3]{index} and \cite[Sect. 4.2]{sectorialextensions} and the references therein. For given $\gamma>0$ we say that a weight function $\omega$ has property $(P_{\omega,\gamma})$ if
\begin{equation*}\label{newindex1}
\exists\;K>1:\;\;\;\limsup_{t\rightarrow+\infty}\frac{\omega(K^{\gamma}t)}{\omega(t)}<K.
\end{equation*}
If $(P_{\omega,\gamma})$ holds for some $K>1$, then $(P_{\omega,\gamma'})$ is satisfied for all $\gamma'\le\gamma$ with the same $K$ since $\omega$ is non-decreasing. Moreover we can restrict to $\gamma>0$, because for $\gamma\le 0$ condition $(P_{\omega,\gamma})$ is satisfied for all $\omega$ (again since $\omega$ is non-decreasing and $K>1$). Then we put
\begin{equation*}\label{newindex2}
\gamma(\omega):=\sup\{\gamma>0: (P_{\omega,\gamma})\;\;\text{is satisfied}\},
\end{equation*}
and if none condition $(P_{\omega,\gamma})$ holds then set $\gamma(\omega):=0$. We recall now some connections between $\gamma(\omega)$ and subadditivity-like conditions for $\omega$.

\begin{remark}\label{secondcomprem}
Let $\omega$ be a weight function.

\begin{itemize}
\item[$(i)$] By \cite[Cor. 2.14]{index} (see also \cite[Lemma 4.2]{sectorialextensions}) we get $\gamma(\omega)>0$ if and only if \hyperlink{om1}{$(\omega_1)$} holds. In particular, for each BMT-weight function $\omega$ we have $\gamma(\omega)>0$.

\item[$(ii)$] By \cite[Thm. 2.11, Cor. 2.13]{index} (see also \cite[Lemma 4.3]{sectorialextensions}) we have that $\gamma(\omega)>1$ if and only if $\omega$ has \hyperlink{omsnq}{$(\omega_{\on{snq}})$}.

\item[$(iii)$] Moreover, by (the proof of) \cite[Thm. 2.11, Cor. 2.13]{index} we know that \hyperlink{omsnq}{$(\omega_{\on{snq}})$} is equivalent to the fact that $\omega\hyperlink{sim}{\sim}\kappa_{\omega}$, where
    \begin{equation}\label{fctkappa}
    \kappa_{\omega}(y):=\int_1^{+\infty}\frac{\omega(yt)}{t^2}dt.
    \end{equation}
    If $\omega$ is also continuous, then by (the proof of) \cite[Prop. 1.3]{MeiseTaylor88} we have that
 $\kappa_{\omega}$ is concave, continuous and $\kappa_{\omega}(0)=\omega(0)\ge 0$; i.e. $\kappa_{\omega}(0)=0$ if and only if $\omega(0)=0$. Hence it follows that also \hyperlink{sub}{$(\omega_{\on{sub}})$} holds for $\kappa_{\omega}$; see e.g. \cite[Lemma 3.8.1 $(1)$]{dissertation}.

 (If continuity for $\omega$ fails, then the proof of \cite[Thm. 2.11 $(i)\Rightarrow(ii)$]{index} yields that we get the same properties for $y\mapsto\widetilde{\kappa}_{\omega}(y):=\int_1^{+\infty}\frac{\kappa_{\omega}(yt)}{t^2}dt$, and so $\omega$ in \eqref{fctkappa} is replaced by $\kappa_{\omega}$. Consequently, in any case one has that $\omega$ is equivalent to a weight satisfying \hyperlink{sub}{$(\omega_{\on{sub}})$}.)

\item[$(iv)$] By definition of $\gamma(\omega)$ we immediately get (see e.g. \cite[$(4.3)$]{sectorialextensions})
\begin{equation}\label{powersubform}
\forall\;s>0:\;\;\;\gamma(\omega^{1/s})=s\gamma(\omega),
\end{equation}
with $\omega^{1/s}(t):=\omega(t^{1/s})$, $t\ge 0$ and $s>0$. Note that $\omega^{1/s}$ is again a weight function; if $\omega$ is a (normalized) BMT-weight function, then each $\omega^{1/s}$, too.

\item[$(v)$] The value of this growth index is preserved under equivalence of weight functions; i.e. if $\sigma,\tau$ are weights such that $\sigma\hyperlink{sim}{\sim}\tau$, then $\gamma(\sigma)=\gamma(\tau)$: This follows from the characterization \cite[Thm. 2.11]{index}, see also \cite[Rem. 2.12]{index}.
\end{itemize}
\end{remark}

Next we summarize some relations between $\gamma(\omega)$ and the notion of (non-)quasianalyticity.

\begin{remark}\label{quasianalyticrem}
We have the following:

\begin{itemize}
\item[$(i)$] Recall that by \cite[Lemma 2.10, Rem. 2.15 $(i)\Rightarrow(iv)$]{index} we know that each quasianalytic weight function $\omega$ has to satisfy $1\le\frac{1}{\gamma(\omega)}(=:\alpha(\omega))$; i.e. $\gamma(\omega)\le 1$.

\item[$(ii)$] Any non-quasianalytic weight function $\omega$ satisfies \hyperlink{om5}{$(\omega_5)$} (and hence \hyperlink{om2}{$(\omega_2)$}) since $\int_a^{+\infty}\frac{\omega(t)}{t^2}dt\ge\omega(a)\int_a^{+\infty}\frac{1}{t^2}dt=\frac{\omega(a)}{a}$ for all $a\ge 1$.
\end{itemize}
\end{remark}

\subsection{Weight sequences and associated weight functions}\label{assofuncsection}
Let $M\in\RR_{>0}^{\NN}$ be given and consider the corresponding sequence of quotients $\mu=(\mu_j)_{j\in\NN}$ by $\mu_j:=\frac{M_j}{M_{j-1}}$, $j\ge 1$, $\mu_0:=1$. $M$ is called {\itshape normalized} if $1=M_0\le M_1$ and $M$ is called {\itshape log-convex} if
$$\forall\;j\in\NN_{>0}:\;M_j^2\le M_{j-1} M_{j+1},$$
equivalently if $(\mu_j)_j$ is non-decreasing. Consider the following set of sequences
$$\hypertarget{LCset}{\mathcal{LC}}:=\{M\in\RR_{>0}^{\NN}:\;M\;\text{is normalized, log-convex},\;\lim_{j\rightarrow+\infty}(M_j)^{1/j}=+\infty\}.$$
We see that $M\in\hyperlink{LCset}{\mathcal{LC}}$ if and only if $1=\mu_0\le\mu_1\le\dots$, $\lim_{j\rightarrow+\infty}\mu_j=+\infty$ (see e.g. \cite[p. 104]{compositionpaper}) and there is a one-to-one correspondence between $M$ and $\mu=(\mu_j)_j$ by taking $M_j:=\prod_{k=0}^j\mu_k$.\vspace{6pt}

$M$ has condition {\itshape moderate growth,} denoted by \hypertarget{mg}{$(\text{mg})$}, if
$$\exists\;C\ge 1\;\forall\;j,k\in\NN:\;M_{j+k}\le C^{j+k+1} M_j M_k.$$
In \cite{Komatsu73} this is denoted by $(M.2)$ and also known under the name {\itshape stability under ultradifferential operators.} $M$ is called {\itshape non-quasianalytic,} denoted by \hypertarget{nq}{$(\text{nq})$} and by condition $(M.3)'$ in \cite{Komatsu73}, if
$$\sum_{j\ge 1}\frac{1}{\mu_j}<+\infty.$$
We call two sequences $M,N\in\RR_{>0}^{\NN}$ {\itshape equivalent} if
$$0<\inf_{j\in\NN_{>0}}\left(\frac{M_j}{N_j}\right)^{1/j}\le\sup_{j\in\NN_{>0}}\left(\frac{M_j}{N_j}\right)^{1/j}<+\infty.$$

Let $M\in\RR_{>0}^{\NN}$ (with $M_0=1$), then the {\itshape associated (weight) function} $\omega_M: \RR\rightarrow\RR\cup\{+\infty\}$ is defined by
\begin{equation*}\label{assofunc}
\omega_M(t):=\sup_{j\in\NN}\log\left(\frac{|t|^j}{M_j}\right)\;\;\;\text{for}\;t\in\RR,\;t\neq 0,\hspace{30pt}\omega_M(0):=0.
\end{equation*}
For an abstract introduction of the associated function we refer to \cite[Chapitre I]{mandelbrojtbook}, see also \cite[Def. 3.1]{Komatsu73} and the more recent work \cite{regularnew}.

\begin{remark}\label{assoweightrem}
Let $M\in\hyperlink{LCset}{\mathcal{LC}}$ be given, then we summarize:
\begin{itemize}
\item[$(*)$] By \cite[Lemma 2.8]{testfunctioncharacterization} and \cite[Lemma 2.4]{sectorialextensions}, see also the references mentioned in the proofs there, $\omega_M$ satisfies all requirements to be a normalized BMT-weight function except \hyperlink{om1}{$(\omega_1)$}.

\item[$(*)$] By \cite[Thm. 3.1]{subaddlike} $\omega_M$ is a normalized BMT-weight function, and so $\omega_M$ has {\itshape in addition} \hyperlink{om1}{$(\omega_1)$}, if and only if $M$ satisfies
    $$\exists\;L\in\NN_{>0}:\;\;\;\liminf_{j\rightarrow+\infty}\frac{(M_{Lj})^{\frac{1}{Lj}}}{(M_j)^{\frac{1}{j}}}>1.$$
    Moreover, by \cite[Prop. 3.4]{subaddlike} it follows that if $M$ has in addition \hyperlink{mg}{$(\on{mg})$}, then $\omega_M$ has \hyperlink{om1}{$(\omega_1)$} if and only if
    \begin{equation}\label{beta3}
    \exists\;Q\in\NN_{\ge 2}:\;\;\;\liminf_{j\rightarrow+\infty}\frac{\mu_{Q_j}}{\mu_j}>1.
    \end{equation}
    \eqref{beta3} appeared crucially in \cite{BonetMeiseMelikhov07} and there it also has been shown that \eqref{beta3} implies \hyperlink{om1}{$(\omega_1)$} for $\omega_M$ even without having \hyperlink{mg}{$(\on{mg})$}; see \cite[Lemma 12, $(2)\Rightarrow(4)$]{BonetMeiseMelikhov07}.

    Note that in \cite[Sect. 2.6 $(2)$]{Meise85}, see also \cite[Rem. 8.9]{BraunMeiseTaylor90}, it has been mentioned that H.-J. Petzsche has already been able to prove the statements \cite[Thm. 3.1, Prop. 3.4]{subaddlike} in a private communication but it seems that a proof never has been published before.

\item[$(*)$] By \cite[Lemma 4.1]{Komatsu73} we have that $\omega_M$ satisfies \hyperlink{omnq}{$(\omega_{\on{nq}})$} if and only if $M$ has \hyperlink{nq}{$(\on{nq})$}.

\item[$(*)$] Finally we mention, see e.g. \cite[Lemma 2.2]{whitneyextensionweightmatrix} and the citations there, that \hyperlink{mg}{$(\on{mg})$} is equivalent to have $\sup_{j\in\NN}\frac{\mu_{2j}}{\mu_j}<+\infty$.
\end{itemize}
\end{remark}

\section{Main results}\label{Abaninmainressect}
We formulate and prove now the main results of this article. In order to proceed we crucially involve the information about the growth index $\gamma(\omega)$. Section \ref{Abaninmaingeneralresultsect} deals with an abstract statement and in Sections \ref{nonquasianalyticweightsequsect} and \ref{quasianalyticweightsequsect} we are investigating in more detail the (non-)quasianalytic weight sequence case.

\subsection{General result for weight functions}\label{Abaninmaingeneralresultsect}
The next technical and crucial result illustrates the intimate connection between the notion of almost subadditivity and the value of $\gamma(\omega)$. Indeed, in this context we study a more general growth relation.

\begin{lemma}\label{crucialindexlemma}
Let $\omega$ be a weight function.
\begin{itemize}
\item[$(i)$] Assume that $\omega$ satisfies
\begin{equation}\label{almostsubaddgeneral}
\exists\;q_0\ge\frac{1}{2}\;\forall\;q>q_0\;\exists\;C_q\ge 1\;\forall\;s,t\ge 0:\;\;\;\omega(s+t)\le q(\omega(s)+\omega(t))+C_q.
\end{equation}
If $q_0>\frac{1}{2}$, then $\gamma(\omega)\ge\frac{\log(2)}{\log(2)+\log(q_0)}=\frac{1}{1+\log(q_0)/\log(2)}>0$ and if $q_0=\frac{1}{2}$, then $\gamma(\omega)=+\infty$.

\item[$(ii)$] If $\omega$ is almost subadditive, i.e. if \eqref{almostsubadd} holds, then $\gamma(\omega)\ge 1$.

\item[$(iii)$] Each almost subadditive and quasianalytic weight function $\omega$ has to satisfy $\gamma(\omega)=1$.
\end{itemize}
\end{lemma}

Concerning $(i)$ note that $\frac{1}{1+\log(q_0)/\log(2)}>0$ because $1+\log(q_0)/\log(2)>0\Leftrightarrow\log(q_0)>-\log(2)\Leftrightarrow q_0>\frac{1}{2}$. Thus, in particular, \eqref{almostsubaddgeneral} yields $\gamma(\omega)>0$ and hence \hyperlink{om1}{$(\omega_1)$}; see $(i)$ in Remark \ref{secondcomprem} but also the proof below.

\demo{Proof}
$(i)$ First, \eqref{almostsubaddgeneral} yields
$$\exists\;q_0\ge\frac{1}{2}\;\forall\;q>q_0\;\exists\;C_q\ge 1\;\forall\;t\ge 0:\;\;\;\omega(2t)\le 2q\omega(t)+C_q,$$
which precisely means
\begin{equation}\label{almostsubadd1}
\exists\;q_0\ge\frac{1}{2}\;\forall\;q>q_0:\;\;\;\limsup_{t\rightarrow+\infty}\frac{\omega(2t)}{\omega(t)}\le 2q;
\end{equation}
in particular \hyperlink{om1}{$(\omega_1)$}. Assume $q_0>\frac{1}{2}$. Let $\epsilon>0$ be small and set $\gamma:=\frac{\log(2)}{\log(2q_0+\epsilon)}$, $K:=2q_0+\epsilon>1$ and so $\gamma=\frac{\log(2)}{\log(K)}$, equivalently $2=K^{\gamma}$. Then, when taking $q$ with $q_0<q<q_0+\frac{\epsilon}{2}$ in \eqref{almostsubadd1}, and so $2q<2q_0+\epsilon=K$, this property verifies $(P_{\omega,\gamma})$ with the aforementioned choice of $K$. One has $0<\gamma<\frac{\log(2)}{\log(2q_0)}$ and when we let $\epsilon\rightarrow 0$ then we get $\gamma\rightarrow\frac{\log(2)}{\log(2q_0)}$ (and $q\rightarrow q_0$). Hence $(P_{\omega,\gamma})$ holds for all $0<\gamma<\frac{\log(2)}{\log(2)+\log(q_0)}$ and by definition this gives $\gamma(\omega)\ge\frac{1}{1+\log(q_0)/\log(2)}$.

However, in general it is not clear that \eqref{almostsubadd1} implies $(P_{\omega,\gamma_0})$ with $\gamma_0:=\frac{1}{1+\log(q_0)/\log(2)}$ and also $\gamma(\omega)>\frac{1}{1+\log(q_0)/\log(2)}$ might be violated.

If $q_0=\frac{1}{2}$, then take $\gamma:=\frac{\log(2)}{\log(1+\epsilon)}$, $K:=1+\epsilon$ and when $\epsilon\rightarrow 0$ by the arguments before we have $(P_{\omega,\gamma})$ for any $\gamma>0$ and so $\gamma(\omega)=+\infty$.\vspace{6pt}

$(ii)$ \eqref{almostsubadd} is precisely \eqref{almostsubaddgeneral} with $q_0=1$.\vspace{6pt}

$(iii)$ If $\omega$ is quasianalytic, then the previous part $(ii)$ and $(i)$ in Remark \ref{quasianalyticrem} yield the conclusion.
\qed\enddemo

\begin{remark}\label{crucialindexlemmarem}
We comment on Lemma \ref{crucialindexlemma}:
\begin{itemize}
\item[$(*)$] The value $\frac{1}{2}$ for the choice of $q_0$ in \eqref{almostsubaddgeneral} is critical: If $0<q_0<\frac{1}{2}$ is allowed, then take $q$ with $q_0<q<\frac{1}{2}$ and we get (see \eqref{almostsubadd1}) that
$$\limsup_{t\rightarrow+\infty}\frac{\omega(2t)}{\omega(t)}\le 2q<1.$$
But this clearly contradicts the standard assumption that $\omega$ is non-decreasing.

\item[$(*)$] As $q_0\rightarrow+\infty$ the lower bound $\frac{1}{1+\log(q_0)/\log(2)}$ tends to $0$.

\item[$(*)$] We give now a partial converse to $(i)$ in Lemma \ref{crucialindexlemma}:

If $\gamma(\omega)>0$, equivalently \hyperlink{om1}{$(\omega_1)$} (see $(i)$ in Remark \ref{secondcomprem}), then since $\omega$ is non-decreasing we estimate for all $s,t\ge 0$ as follows
\begin{equation}\label{subbaddvsomega}
\omega(s+t)\le\omega(2\max\{s,t\})\le L\omega(\max\{s,t\})+L\le L(\omega(t)+\omega(s))+L.
\end{equation}
Consequently, $\gamma(\omega)>0$ yields \eqref{almostsubaddgeneral} for $q_0=L$, with $L$ being the constant such that $\omega(2t)\le L\omega(t)+L$ for all $t\ge 0$; see \hyperlink{om1}{$(\omega_1)$}. Note that $L\ge 1$ because $L<1$ yields a contradiction to the fact that $\omega$ is a weight function (so $\lim_{t\rightarrow+\infty}\omega(t)=+\infty$ and $\omega$ is non-decreasing).

Similarly, it is not clear that \eqref{almostsubadd1} implies \eqref{almostsubaddgeneral}: The analogous estimate as in \eqref{subbaddvsomega} gives that \eqref{almostsubadd1} for $q_0$ yields \eqref{almostsubaddgeneral} for $2q_0$.
\end{itemize}
\end{remark}

Using this information we are able to show the following main statement.

\begin{theorem}\label{Abaninmainthm}
For any weight function $\omega$ with $\gamma(\omega)<+\infty$ we get:
\begin{itemize}
\item[$(*)$] There exist uncountable infinite many pair-wise different (i.e. non-equivalent) weight functions $\tau_a:=\omega^{1/a}$, $a\in(0,\gamma(\omega)^{-1})=:\mathcal{I}_{\omega}$ (with $\mathcal{I}_{\omega}=(0,+\infty)$ if $\gamma(\omega)=0$), such that each $\tau_a$ has the following property:

There does not exist a weight function $\sigma$ being equivalent to $\tau_a$ and such that $\sigma$ is almost subadditive (see \eqref{almostsubadd}). In particular, each $\tau_a$ itself is violating \eqref{almostsubadd}.

\item[$(*)$] If $\omega$ is a (normalized) BMT-weight function, then each $\tau_a$, too.
\end{itemize}
\end{theorem}

Roughly speaking this result means that ``many'' BMT-weight functions are not equivalent to an almost subadditive function in the sense of \cite{Abanin10}; to each weight with finite $\gamma(\omega)$ we can assign an uncountable infinite family of pair-wise different equivalence classes of weights such that in each of these classes no almost subadditive weight exists.

\demo{Proof}
Take any $a>0$ with $\gamma(\omega)<\frac{1}{a}$, then set $\tau_a:=\omega^{1/a}$ and by \eqref{powersubform} we get
\begin{equation}\label{Abaninmainthmequ}
\gamma(\tau_a)=\gamma(\omega^{1/a})=a\gamma(\omega)<1.
\end{equation}
Note: If $\gamma(\omega)=0$, then any choice for $a>0$ in \eqref{Abaninmainthmequ} is admissible and $\gamma(\tau_a)=0$ for all $a>0$. This verifies the interpretation $\mathcal{I}_{\omega}=(0,+\infty)$ if $\gamma(\omega)=0$. (On the other hand, in order to choose $a$ the assumption $\gamma(\omega)<+\infty$ is indispensable. In the case $\gamma(\omega)=+\infty$, which is excluded by assumption, we put $\mathcal{I}_{\omega}:=\emptyset$.)

Consider now for such $a>0$ (fixed) the function $\tau_a=\omega^{1/a}$. By $(ii)$ in Lemma \ref{crucialindexlemma} and \eqref{Abaninmainthmequ} we have that $\omega_a$ is violating \eqref{almostsubadd}. Since $\gamma(\cdot)$ is preserved under equivalence of weight functions, see $(v)$ in Remark \ref{secondcomprem}, we have that each weight $\sigma$ with $\sigma\hyperlink{sim}{\sim}\tau_a$ satisfies $\gamma(\sigma)<1$ and hence also cannot have \eqref{almostsubadd}.

Recall that if $\omega$ is a (normalized) BMT-weight, then each $\tau_a=\omega^{1/a}$, too.\vspace{6pt}

Finally, let $a_1>a_2>0$ satisfy $\gamma(\omega)<\frac{1}{a_1}<\frac{1}{a_2}$ and let these choices be arbitrary but from now on fixed. We prove that $\omega^{1/a_1}$ and $\omega^{1/a_2}$ are not equivalent:

On the one hand, since $\omega$ is assumed to be non-decreasing, we clearly have $\omega^{1/a_1}(t)\le\omega^{1/a_2}(t)\Leftrightarrow\omega(t^{1/a_1})\le\omega(t^{1/a_2})$ for any $0<a_2<a_1$ and all $t\ge 1$. But, on the other hand, if the weights are equivalent, i.e. if $\omega^{1/a_1}\hyperlink{sim}{\sim}\omega^{1/a_2}$, then
$$\exists\;C\ge 1\;\forall\;t\ge 0:\;\;\;\omega(t^{1/a_2})\le C\omega(t^{1/a_1})+C,$$
and so
\begin{equation}\label{om7var}
\exists\;C\ge 1\;\forall\;t\ge 0:\;\;\;\omega(t^{a_1/a_2})\le C\omega(t)+C.
\end{equation}
We show that this is equivalent to
\begin{equation}\label{om7}
\exists\;C\ge 1\;\forall\;t\ge 0:\;\;\;\omega(t^2)\le C\omega(t)+C.
\end{equation}
\eqref{om7var}$\Rightarrow$\eqref{om7}: If $a_1\ge 2a_2$, then the claim follows because $\omega$ is non-decreasing. If $2a_2>a_1>a_2$, then an $n$-times iteration with $n\in\NN_{>0}$ (minimal) such that $(a_1/a_2)^n\ge 2$ yields the conclusion.

\eqref{om7}$\Rightarrow$\eqref{om7var}: Similarly, if $2a_2\ge a_1>a_2$, then the implication is clear. If $a_1>2a_2$, then an $n$-times iteration with $n\in\NN_{>0}$ (minimal) such that $a_1/a_2\le 2^n$ yields the conclusion.\vspace{6pt}

Finally, we verify that \eqref{om7} implies $\gamma(\omega)=+\infty$ for any weight function: Fix some $K>1$ arbitrary and take any $\gamma>0$. (As we are going to see the choice for $K$ can be taken uniformly for all $\gamma>0$.) Then we can find $t_{\gamma,K}:=K^{\gamma}>0$ such that $\omega(K^{\gamma}t)\le\omega(t^2)$ for all $t\ge t_{\gamma,K}$ since $\omega$ is non-decreasing and so $\omega(K^{\gamma}t)\le C\omega(t)+C$ for all $t\ge t_{\gamma,K}$. Hence it suffices to choose $K:=C-\epsilon$ with $C$ denoting the constant from \eqref{om7} and $\epsilon>0$ small and fixed.

But this is a contradiction to the basic assumption $\gamma(\omega)<+\infty$. (This step of the proof should be compared with \cite[$(A.1)$, p. 2171]{sectorialextensions1} and the proof of \cite[Lemma A.1, p. 2173]{sectorialextensions1}; formally there we have dealt with more specific weight functions but the additional assumptions on the weights are superfluous.)
\qed\enddemo

\begin{example}\label{promienentexample}
Consider the Gevrey weights $\omega_s(t):=t^{1/s}$, $s>0$. Each $\omega_s$ is a BMT-weight function; if $s>1$ then $\omega_s$ is non-quasianalytic and for all $0<s\le 1$ the weight is quasianalytic. It is straight-forward to see that $\gamma(\omega_s)=s$. Fix $s>0$ and by Theorem \ref{Abaninmainthm} for each $a\in(0,s^{-1})$, and so $as<1,$ the weight $\tau_a=\omega_s(t^{1/a})=t^{1/(sa)}$ has the desired property.

Note that for each $s\ge 1$ we have that \hyperlink{sub}{$(\omega_{\on{sub}})$} is valid and for each $s>1$ even \hyperlink{omsnq}{$(\omega_{\on{snq}})$}; see $(iii)$ in Remark \ref{secondcomprem}.

Summarizing, for any $s\ge 1$ the weight $\omega_s$ is (almost) subadditive but, if $0<s<1$, then there does not exist a weight function $\tau$ such that $\tau$ is almost subadditive and $\tau\hyperlink{sim}{\sim}\omega_s$. Compare this also with $(iii)$ in Lemma \ref{crucialindexlemma}. However, since for each such small $s$ the weight is quasianalytic, formally this example does not solve the question from \cite[Rem. 1.5]{Abanin10}.\vspace{6pt}

$\sigma_s(t):=\max\{0,(\log(t))^s\}$, $s>1$, is another prominent standard example for a (normalized) BMT-weight function. However, here $\gamma(\sigma_s)=+\infty$ for any $s>1$ and hence Theorem \ref{Abaninmainthm} cannot be applied.
\end{example}

\begin{remark}\label{Abaninmainthmrem}
The (standard) assumptions on the weight $\omega$ in Theorem \ref{Abaninmainthm} are more general than in \cite{Abanin10} and, in particular, this statement applies to quasianalytic weights as well since by $(i)$ in Remark \ref{quasianalyticrem} in this case $\gamma(\omega)\le 1$ is valid.

Moreover, if instead of almost subadditivity we consider the more general condition \eqref{almostsubaddgeneral}, then Theorem \ref{Abaninmainthm} transfers to this situation with the analogous proof; only the interval $\mathcal{I}_{\omega}$ changes. If $q_0>\frac{1}{2}$, then we get $\mathcal{I}_{\omega}=(0,\frac{1}{1+\log(q_0)/\log(2)}\gamma(\omega)^{-1})$ (again with $\mathcal{I}_{\omega}=(0,+\infty)$ provided that $\gamma(\omega)=0$): choose $a>0$ with $\left(1+\frac{\log(q_0)}{\log(2)}\right)\gamma(\omega)<\frac{1}{a}$ and so by \eqref{powersubform}
$$\gamma(\tau_a)=\gamma(\omega^{1/a})=a\gamma(\omega)<\frac{\log(2)}{\log(2)+\log(q_0)}.$$
The rest follows as before. If $q_0=\frac{1}{2}$, then we have $\mathcal{I}_{\omega}=(0,+\infty)$ and any choice for $a>0$ is valid.
\end{remark}

\subsection{The non-quasianalytic weight sequence case}\label{nonquasianalyticweightsequsect}
From now on we focus on the case $\omega=\omega_M$ for a specific sequence $M\in\RR_{>0}^{\NN}$. The proofs are involved and use technical preparations and therefore we split it into several steps. However, the ideas and constructions can be useful for further applications.

\begin{theorem}\label{Abaninexamplethm}
There exist uncountable infinite many pair-wise different (i.e. non-equivalent) weight sequences $M\in\hyperlink{LCset}{\mathcal{LC}}$ and normalized and non-quasianalytic BMT-weight functions $\omega=\omega_M$ such that each of these functions has the following property:

There does not exist a weight function $\tau$ such that $\tau$ is equivalent to $\omega_M$ and $\tau$ is almost subadditive.
\end{theorem}

Note that $\omega_s$ in Example \ref{promienentexample} corresponds, up to equivalence, to $\omega_{G^s}$ with $G^s:=(j!^s)_{j\in\NN}$. But for $0<s\le 1$ the weights $\omega_s$, $\omega_{G^s}$ (and the sequences $G^s$) are quasianalytic and hence not an appropriate choice.

\demo{Proof}
{\itshape Step I} - Comments on sufficient growth and regularity properties.\vspace{6pt}

By Lemma \ref{crucialindexlemma} in order to conclude it suffices to construct a normalized BMT-weight function $\omega$ which is
\begin{itemize}
\item[$(*)$] non-quasianalytic and

\item[$(*)$] such that $0<\gamma(\omega)<1$.
\end{itemize}

In order to achieve this goal the idea is now to consider $\omega\equiv\omega_M$ with $M\in\hyperlink{LCset}{\mathcal{LC}}$ satisfying explicit growth and regularity requirements and to define $M$ explicitly.\vspace{6pt}

{\itshape Step II} - Sufficient growth and regularity properties for the weight sequence.\vspace{6pt}

We want to find $M\in\hyperlink{LCset}{\mathcal{LC}}$ such that
\begin{itemize}
\item[$(a)$] $M$ satisfies \eqref{beta3},

\item[$(b)$] $M$ satisfies \hyperlink{nq}{$(\on{nq})$},

\item[$(c)$] $M$ satisfies \hyperlink{mg}{$(\on{mg})$},

\item[$(d)$] $M$ satisfies
\begin{equation}\label{gammaMrelation}
\exists\;0<\beta<1\;\forall\;\ell\in\NN_{\ge 2}:\;\;\;\liminf_{j\rightarrow+\infty}\frac{\mu_{\ell j}}{\mu_j}\le\ell^{\beta}.
\end{equation}
\end{itemize}
Because for such $M$ we get the following properties; for the precise definition of the growth index $\gamma(M)$, its comparison with $\gamma(\omega_M)$ and more details we refer to \cite[Sect. 3-5]{index}:

\begin{itemize}
	\item[$(*)$] By $(a)$ and $(b)$ the function $\omega_{M}$ is a non-quasianalytic (and normalized) BMT-weight (see Remark \ref{assoweightrem}).
	
	\item[$(*)$] By $(c)$ and \cite[Cor. 4.6 $(iii)$]{index} we get $\gamma(M)=\gamma(\omega_M)$. Note: Without having \hyperlink{mg}{$(\on{mg})$} we only know $\gamma(M)\le\gamma(\omega_M)$ and the difference can be large in general; see \cite[Cor. 4.6 $(i)$, Sect. 5]{index}.
	
	\item[$(*)$] By $(d)$ and \cite[Thm. 3.11 $(v)\Leftrightarrow(viii)$]{index} we see that $\gamma(M)\le\beta<1$. Recall that $M\in\hyperlink{LCset}{\mathcal{LC}}$ implies $\gamma(M)\ge 0$.

\item[$(*)$] Note that in \cite{index} a weight sequence denotes a sequence satisfying all requirements from the class \hyperlink{LCset}{$\mathcal{LC}$} except necessarily $M_0\le M_1$ (see \cite[Sect. 3.1]{index}). Moreover, concerning the quotients an index-shift appears; i.e. $m_p$ in \cite{index} is corresponding to $\mu_{p+1}$. However, this does not effect the value of the index and hence the conclusion; see \cite[Rem. 3.8]{index} and the comments after \cite[Cor. 3.12]{index}.
\end{itemize}

Altogether, $\omega_{M}$ is a (normalized) non-quasianalytic BMT-weight function satisfying $\gamma(\omega_{M})=\gamma(M)<1$ and by {\itshape Step I} this function provides a convenient example. Thus it remains to construct explicitly $M\in\hyperlink{LCset}{\mathcal{LC}}$ such that $(a)-(d)$ holds.\vspace{6pt}

{\itshape Step III} - Definition of the sequence $M$.\vspace{6pt}

We introduce $M$ in terms of its quotient sequence $\mu=(\mu_k)_{k\in\NN}$ (with $\mu_0:=1$) and so set $M_j:=\prod_{k=0}^j\mu_k$, $j\in\NN$. For this, first fix an arbitrary real constant $A>2$. Second, we consider three auxiliary sequences (of positive integers) $(a_j)_{j\in\NN_{>0}}$, $(b_j)_{j\in\NN_{>0}}$ and $(d_j)_{j\in\NN_{>0}}$ with $a_j<b_j<a_{j+1}$ for all $j\in\NN_{>0}$ and being defined as follows:
\begin{equation}\label{AbaninStepIVequ}
a_1:=1,\hspace{20pt}b_j:=2^{d_j}a_j,\hspace{20pt}a_{j+1}:=2^jb_j,\hspace{20pt}j\in\NN_{>0}.
\end{equation}
For the sequence $(d_j)_{j\in\NN_{>0}}$ we assume that it is strictly increasing with $d_1\ge 1$ and satisfying the following growth restriction:
\begin{equation}\label{AbaninStepIVequ5}
\forall\;j\in\NN_{>0}:\;\;\;\left(\frac{2}{A}\right)^{d_j}\le\frac{1}{2^{\frac{j}{2}+1}}.
\end{equation}
This choice for $d_j$ is possible since $A>2$.

For notational reasons we split each interval $(a_j,b_j]$, $j\in\NN_{>0}$, into $d_j$-many subintervals $I_{a_j,i}:=(2^ia_j,2^{i+1}a_j]$, $i\in\NN$ and $0\le i\le d_j-1$. Similarly, each interval $(b_j,a_{j+1}]$, $j\in\NN_{>0}$, is split into $j$-many subintervals $I_{b_j,i}:=(2^ib_j,2^{i+1}b_j]$, $i\in\NN$ and $0\le i\le j-1$. To be formally precise, we only consider all integers belonging to these intervals.

We define $\mu$ on $a_1$ and on the intervals $(a_j,b_j]$, $j\in\NN_{>0}$, as follows:
\begin{equation}\label{AbaninStepIVequ1}
\mu_{a_1}:=2(=2a_1),\hspace{30pt}\mu_k:=A\mu_{2^ia_j}=A^{i+1}\mu_{a_j},\;\;\;k\in I_{a_j,i},\;\;\;0\le i\le d_j-1,\;j\in\NN_{>0}.
\end{equation}
This means that iteratively on $I_{a_j,i}$, $i\ge 1$, the sequence of quotients is given by multiplying the already defined quotients on $I_{a_j,i-1}$ with the constant $A$ and $\mu$ is constant on each $I_{a_j,i}$. And the values on $I_{a_j,0}$ are obtained by multiplying $\mu_{a_j}$ with $A$. In particular, we get
\begin{equation}\label{AbaninStepIVequ1impl}
\forall\;j\in\NN_{>0}:\;\;\;\mu_{b_j}=A^{d_j}\mu_{a_j}.
\end{equation}

Next, on each $(b_j,a_{j+1}]$, $j\in\NN_{>0}$, we introduce $\mu$ as follows:
\begin{equation}\label{AbaninStepIVequ4}
\mu_k:=\sqrt{2}^{\frac{1}{2^ib_j}}\mu_{k-1},\;\;\;2^ib_j\le k-1<k\le 2^{i+1}b_j,\hspace{20pt}0\le i\le j-1,\;j\in\NN_{>0}.
\end{equation}
Note: When considering in \eqref{AbaninStepIVequ4} the choices $i=j-1$ and $k=2^{i+1}b_j=2^jb_j=a_{j+1}$ this gives the value of $\mu_{a_{j+1}}$ and so the definition of $M$ is complete.

Since $\frac{\mu_{2k}}{\mu_k}=\frac{\mu_{2k}}{\mu_{2k-1}}\frac{\mu_{2k-1}}{\mu_{2k-2}}\cdots\frac{\mu_{k+1}}{\mu_k}$ for all $k\in\NN_{>0}$ gives a product of $k$-many factors, \eqref{AbaninStepIVequ4} yields
\begin{equation}\label{AbaninStepIVequ4impl}
\forall\;j\in\NN_{>0}\;\forall\;0\le i\le j-1:\;\;\;\frac{\mu_{2^{i+1}b_j}}{\mu_{2^ib_j}}=\sqrt{2},
\end{equation}
and which implies by iteration
\begin{equation}\label{AbaninStepIVequ4imp2}
\forall\;j\in\NN_{>0}\;\forall\;0\le i\le j-1:\;\;\;\frac{\mu_{2^{i+1}b_j}}{\mu_{b_j}}=2^{\frac{i+1}{2}}.
\end{equation}
In particular, the choice $i=j-1$ gives
\begin{equation}\label{AbaninStepIVequ4imp3}
\forall\;j\in\NN_{>0}:\;\;\;\mu_{a_{j+1}}=2^{\frac{j}{2}}\mu_{b_j}.
\end{equation}
Let us now show by induction the following consequence of \eqref{AbaninStepIVequ5} which is needed in the estimates below:
\begin{equation}\label{AbaninStepIVequ6}
\forall\;j\in\NN_{>0}:\;\;\;\frac{a_j}{\mu_{a_j}}\le\frac{1}{2^j}.
\end{equation}
The case $j=1$ follows by \eqref{AbaninStepIVequ1}. Moreover, combining \eqref{AbaninStepIVequ}, \eqref{AbaninStepIVequ1impl} \eqref{AbaninStepIVequ4imp3} with the induction hypothesis and finally with \eqref{AbaninStepIVequ5} gives for all $j\in\NN_{>0}$
$$\frac{a_{j+1}}{\mu_{a_{j+1}}}=\frac{2^jb_j}{2^{\frac{j}{2}}\mu_{b_j}}=\frac{2^j2^{d_j}a_j}{2^{\frac{j}{2}}A^{d_j}\mu_{a_j}}=\frac{a_j}{\mu_{a_j}}2^{\frac{j}{2}}\left(\frac{2}{A}\right)^{d_j}\le\frac{1}{2^{\frac{j}{2}}}\left(\frac{2}{A}\right)^{d_j}\le\frac{1}{2^{j+1}},$$
as desired.\vspace{6pt}

{\itshape Step IV} - Verification of the desired properties for $M$. We divide this step into several claims.\vspace{6pt}

{\itshape Claim a:} $M\in\hyperlink{LCset}{\mathcal{LC}}$. This is obvious since $M_0=\mu_0=1$ and $\frac{M_1}{M_0}=\mu_1=\mu_{a_1}=2$ ensure normalization. Clearly $j\mapsto\mu_j$ is by construction non-decreasing and satisfying $\lim_{j\rightarrow+\infty}\mu_j=+\infty$.\vspace{6pt}

{\itshape Claim b:} $M$ satisfies \hyperlink{nq}{$(\on{nq})$}.\vspace{6pt}

On the one hand, by recalling definition \eqref{AbaninStepIVequ1} and by \eqref{AbaninStepIVequ6} we have for all $j\in\NN_{>0}$:
\begin{align*}
\sum_{k=a_j+1}^{b_j}\frac{1}{\mu_k}&=\sum_{i=0}^{d_j-1}\sum_{k\in I_{a_j,i}}\frac{1}{\mu_k}=\sum_{i=0}^{d_j-1}\frac{2^ia_j}{A^{i+1}\mu_{a_j}}=\frac{a_j}{A\mu_{a_j}}\sum_{i=0}^{d_j-1}\left(\frac{2}{A}\right)^i\le\frac{a_j}{A\mu_{a_j}}\sum_{i=0}^{+\infty}\left(\frac{2}{A}\right)^i
\\&
=\frac{a_j}{A\mu_{a_j}}\frac{A}{A-2}\le\frac{1}{2^j}\frac{1}{A-2}.
\end{align*}
On the other hand, by \eqref{AbaninStepIVequ}, \eqref{AbaninStepIVequ4} and by combining \eqref{AbaninStepIVequ5} and \eqref{AbaninStepIVequ6} we get for all $j\in\NN_{>0}$:
$$\sum_{k=b_j+1}^{a_{j+1}}\frac{1}{\mu_k}=\sum_{i=0}^{j-1}\sum_{k\in I_{b_j,i}}\frac{1}{\mu_k}\le\frac{a_{j+1}-b_j}{\mu_{b_j}}=\frac{b_j(2^j-1)}{A^{d_j}\mu_{a_j}}=\left(\frac{2}{A}\right)^{d_j}\frac{a_j(2^j-1)}{\mu_{a_j}}\le\frac{2^j-1}{2^{j+\frac{j}{2}+1}}\le\frac{1}{2^{\frac{j}{2}+1}}.$$
Summarizing, we obtain
\begin{align*}
\sum_{k\ge 1}\frac{1}{\mu_k}&=\sum_{k\ge a_1}\frac{1}{\mu_k}=\frac{1}{2}+\sum_{k>a_1}\frac{1}{\mu_k}\le\frac{1}{2}+\frac{1}{A-2}\sum_{j\ge 1}\frac{1}{2^j}+\frac{1}{2}\sum_{j\ge 1}\frac{1}{2^{\frac{j}{2}}}
\\&
=\frac{A}{2(A-2)}+\frac{1}{2}\frac{1}{\sqrt{2}-1}<+\infty,
\end{align*}
and the claim is shown.\vspace{6pt}

{\itshape Claim c:} $M$ satisfies \hyperlink{mg}{$(\on{mg})$}.\vspace{6pt}

Recall that in order to conclude we have to show $\sup_{k\in\NN}\frac{\mu_{2k}}{\mu_k}<+\infty$ and according to the definition of $M$ we have to distinguish between several cases (also required in {\itshape Claim d} below):\vspace{6pt}

\begin{itemize}
\item[$(i)$] For $k=1$ we get $\frac{\mu_{2k}}{\mu_k}=\frac{\mu_{2a_1}}{\mu_{a_1}}=A$, see \eqref{AbaninStepIVequ} and \eqref{AbaninStepIVequ1}.

\item[$(ii)$] When $k\in I_{a_j,i}$ for $j\in\NN_{>0}$ and $0\le i\le d_j-2$ arbitrary, then we get $2k\in I_{a_j,i+1}$ and so $\frac{\mu_{2k}}{\mu_k}=A$, see again \eqref{AbaninStepIVequ1}.

\item[$(iii)$] When $k\in I_{a_j,d_j-1}$, $j\in\NN_{>0}$ arbitrary, then $2k\in I_{b_j,0}$ and there exists $2\le\ell\le b_j$ such that $b_j+\ell=2k$. This gives $$\frac{\mu_{2k}}{\mu_k}=\frac{\mu_{2k}}{\mu_{b_j}}\frac{\mu_{b_j}}{\mu_k}=\sqrt{2}^{\frac{\ell}{b_j}},$$
    see \eqref{AbaninStepIVequ1} and \eqref{AbaninStepIVequ4}, and this equality implies the estimate $\frac{\mu_{2k}}{\mu_k}\le\sqrt{2}$.

\item[$(iv)$] When $k\in I_{b_j,j-1}$, $j\in\NN_{>0}$ arbitrary, then $2k\in I_{a_{j+1},0}$ and there exists $0\le\ell\le 2^{j-1}b_j-1$ such that $k+\ell=a_{j+1}$. This implies
    $$\frac{\mu_{2k}}{\mu_k}=\frac{\mu_{a_{j+1}}}{\mu_k}\frac{\mu_{2k}}{\mu_{a_{j+1}}}=\sqrt{2}^{\frac{\ell}{2^{j-1}b_j}}A,$$
see \eqref{AbaninStepIVequ1} and \eqref{AbaninStepIVequ4}, and this equality yields the estimate $\frac{\mu_{2k}}{\mu_k}\le\sqrt{2}A$.

\item[$(v)$] Finally, let $k\in I_{b_j,i}$ with $0\le i\le j-2$ and $j\in\NN_{>0}$ arbitrary. In this case $2k\in I_{b_j,i+1}$ and there exists $0\le\ell\le 2^ib_j-1$ such that $k+\ell=2^{i+1}b_j$ and there exists $2\le\ell'\le 2^{i+1}b_j$ such that $2^{i+1}b_j+\ell'=2k$ (and $\ell,\ell'$ satisfy the relation $k+\ell+\ell'=2k$). Thus we get
    $$\frac{\mu_{2k}}{\mu_k}=\frac{\mu_{2^{i+1}b_j}}{\mu_k}\frac{\mu_{2k}}{\mu_{2^{i+1}b_j}}=\sqrt{2}^{\frac{\ell}{2^ib_j}}\sqrt{2}^{\frac{\ell'}{2^{i+1}b_j}},$$
    see \eqref{AbaninStepIVequ4}. This equality implies the estimate $$\frac{\mu_{2k}}{\mu_k}\le\sqrt{2}^{\frac{\ell}{2^ib_j}}\sqrt{2}^{\frac{\ell'}{2^ib_j}}=\sqrt{2}^{\frac{\ell+\ell'}{2^ib_j}}=\sqrt{2}^{\frac{k}{2^ib_j}}\le\sqrt{2}^{\frac{2^{i+1}b_j}{2^ib_j}}=2.$$
\end{itemize}
Summarizing, we have verified $\frac{\mu_{2k}}{\mu_k}\le\sqrt{2}A$ for all $k\in\NN_{>0}$ and we are done. (The case $k=0$ is trivial and gives $\frac{\mu_{2k}}{\mu_k}=1$.)\vspace{6pt}

{\itshape Claim d:} $M$ satisfies \eqref{beta3}.\vspace{6pt}

We verify \eqref{beta3} for the choice $Q=4$ and use the equalities and explanations from $(i)-(v)$ in {\itshape Claim c} above. (In fact even $Q=2$ is sufficient to treat all cases except $(iii)$.)

For all $k$ as considered in $(ii)$ in {\itshape Claim c} we get with $Q=2$ the same equality and for $k=1$ in case $(i)$, too. The equality in $(iv)$ in {\itshape Claim c} yields $\frac{\mu_{2k}}{\mu_k}=\sqrt{2}^{\frac{\ell}{2^{j-1}b_j}}A\ge A$ and analogous estimates as in $(v)$ in {\itshape Claim c} give $$\frac{\mu_{2k}}{\mu_k}\ge\sqrt{2}^{\frac{\ell}{2^{i+1}b_j}}\sqrt{2}^{\frac{\ell'}{2^{i+1}b_j}}=\sqrt{2}^{\frac{\ell+\ell'}{2^{i+1}b_j}}=\sqrt{2}^{\frac{k}{2^{i+1}b_j}}\ge\sqrt{2}^{\frac{2^ib_j}{2^{i+1}b_j}}=2^{1/4}.$$

Finally, let $k$ be as in case $(iii)$ in {\itshape Claim c.} Then $4k\in I_{b_j,1}$ and so
$$\frac{\mu_{4k}}{\mu_k}=\frac{\mu_{4k}}{\mu_{b_j}}\ge\frac{\mu_{2b_j}}{\mu_{b_j}}=\sqrt{2},$$
see \eqref{AbaninStepIVequ4impl}.

Summarizing, even $\inf_{j\rightarrow+\infty}\frac{\mu_{4j}}{\mu_j}\ge 2^{1/4}>1$ is verified and thus the claim is shown.\vspace{6pt}

{\itshape Claim e:} $M$ satisfies \eqref{gammaMrelation}.\vspace{6pt}

Let $\ell\in\NN_{\ge 2}$ be given, arbitrary but from now on fixed. Consider $n\in\NN_{>0}$ such that $2^n\le\ell<2^{n+1}$ is valid and distinguish: When $\ell=2^n$ for some $n\in\NN_{>0}$, in particular when $\ell=2$, then \eqref{AbaninStepIVequ4imp2} implies for all $j\ge n$
$$\frac{\mu_{\ell b_j}}{\mu_{b_j}}=\frac{\mu_{2^nb_j}}{\mu_{b_j}}=2^{n/2}=\sqrt{\ell},$$
which verifies the $\liminf$-condition in \eqref{gammaMrelation} with $\beta=\frac{1}{2}$.

When $2^n<\ell<2^{n+1}$, then for all $j\ge n+1$ property \eqref{AbaninStepIVequ4imp2} (and the fact that $k\mapsto\mu_k$ is non-decreasing) yields
$$\frac{\mu_{\ell b_j}}{\mu_{b_j}}\le\frac{\mu_{2^{n+1}b_j}}{\mu_{b_j}}=2^{\frac{n+1}{2}}=\sqrt{2}2^{\frac{n}{2}}\le\sqrt{2\ell}\le\ell^{\beta},$$
where we have put $\beta:=\frac{1}{2}+\frac{1}{2}\frac{\log(2)}{\log(3)}$ and therefore $\frac{1}{2}<\beta<1$. Note that with this choice for $\beta$ the last estimate above is equivalent to $2\le\ell^{\log(2)/\log(3)}$ and so to $\log(3)\le\log(\ell)$ which is obviously true for all $\ell$ under consideration since we have $\ell\in\NN_{\ge 3}$.

Summarizing, \eqref{gammaMrelation} is verified with $\beta:=\frac{1}{2}+\frac{1}{2}\frac{\log(2)}{\log(3)}$.\vspace{6pt}

{\itshape Step V} - There exist uncountable infinite many pair-wise different weight sequences/functions satisfying the required properties.\vspace{6pt}

Given a parameter $A>2$, then in view of \eqref{AbaninStepIVequ5} for all $A'>A$ we can take the same auxiliary sequence $(d_j)_{j\in\NN_{>0}}$ and hence also the sequences $(a_j)_{j\in\NN_{>0}}$ and $(b_j)_{j\in\NN_{>0}}$ remain unchanged; see \eqref{AbaninStepIVequ}. So, when fixing $A_0>2$ and varying $A$ in $[A_0,+\infty)$, then in the construction all claims stay valid and it is easy to compare the quotient sequences since the positions where the growth behavior of this sequence is changing remain invariant.

Consider now $A'>A>2$ arbitrary but from now on fixed, then we write $M^A$ and $M^{A'}$ for the corresponding sequences and similarly $\mu^A$ and $\mu^{A'}$ for the corresponding sequences of quotients. By \eqref{AbaninStepIVequ1impl} and \eqref{AbaninStepIVequ4imp3} one has
$$\forall\;j\in\NN_{>0}:\;\;\;\frac{\mu^{A'}_{a_{j+1}}}{\mu^A_{a_{j+1}}}=\frac{2^{\frac{j}{2}}(A')^{d_j}\mu^{A'}_{a_j}}{2^{\frac{j}{2}}A^{d_j}\mu^A_{a_j}}=\left(\frac{A'}{A}\right)^{d_j}\frac{\mu^{A'}_{a_j}}{\mu^A_{a_j}}.$$
Iterating this procedure we immediately get
$$\forall\;j\in\NN_{>0}:\;\;\;\frac{\mu^{A'}_{a_{j+1}}}{\mu^A_{a_{j+1}}}=\left(\frac{A'}{A}\right)^{d_j+d_{j-1}+\dots+d_1}\frac{\mu^{A'}_{a_1}}{\mu^A_{a_1}}=\left(\frac{A'}{A}\right)^{d_j+d_{j-1}+\dots+d_1},$$
hence $\limsup_{k\rightarrow+\infty}\frac{\mu^{A'}_k}{\mu^A_k}=+\infty$ is verified.

Next, we claim that this property implies that $M^A$ and $M^{A'}$ are not equivalent: Both sequences satisfy \hyperlink{mg}{$(\on{mg})$} (see {\itshape Claim c}) and since $M^{A'}\in\hyperlink{LCset}{\mathcal{LC}}$ condition \hyperlink{mg}{$(\on{mg})$} is equivalent to $\sup_{k\in\NN_{>0}}\frac{\mu^{A'}_k}{(M^{A'}_k)^{1/k}}<+\infty$, see e.g. again \cite[Lemma 2.2]{whitneyextensionweightmatrix}. Moreover, by log-convexity and normalization we clearly have $(M^A_k)^{1/k}\le\mu^A_k$ for all $k\in\NN_{>0}$ and hence
\begin{equation}\label{nonequivequ}
\sup_{j\in\NN_{>0}}\left(\frac{M^{A'}_j}{M^A_j}\right)^{1/j}=+\infty.
\end{equation}
Thus the sequences are not equivalent.

Now recall that $M\in\hyperlink{LCset}{\mathcal{LC}}$ has \hyperlink{mg}{$(\on{mg})$} if and only if $M$ is equivalent to some/all auxiliary sequences $\widetilde{M}^c$ with $\widetilde{M}^c_j:=(M_{cj})^{1/c}$, $c\in\NN_{>0}$; see \cite[Lemma 2.2]{subaddlike}.

Consequently, since $M^A$ has \hyperlink{mg}{$(\on{mg})$}, we see that the estimate
\begin{equation}\label{64equ}
\exists\;c\in\NN_{>0}\;\exists\;B\ge 1\;\forall\;j\in\NN:\;\;\;M^{A'}_j\le B(M^A_{cj})^{1/c}
\end{equation}
has to be violated: Otherwise $M^{A'}_j\le B(M^A_{cj})^{1/c}\le BC^jM^A_j$ for some $B,C\ge 1$ and all $j\in\NN$ and this contradicts \eqref{nonequivequ}. \eqref{64equ} is precisely \cite[$(6.4)$]{PTTvsmatrix} for $N\equiv M^{A'}$ and $M\equiv M^A$. Therefore, the first part in \cite[Lemma 6.5]{PTTvsmatrix} implies that $\omega_{M^A}(t)=O(\omega_{M^{A'}}(t))$ as $t\rightarrow+\infty$ fails and so $\omega_{M^{A}}$ and $\omega_{M^{A'}}$ are not equivalent.
\qed\enddemo

\subsection{The quasianalytic weight sequence case}\label{quasianalyticweightsequsect}
The aim is to transfer Theorem \ref{Abaninexamplethm} also to quasianalytic BMT-weight functions; so we illustrate independently that the question whether there does always exist an equivalent almost subadditive BMT-weight function is not depending on the notion of (non-)quasianalyticity.

First we comment in more detail on the growth of $j\mapsto\frac{\mu_j}{j}$ in the example of Theorem \ref{Abaninexamplethm}. Given $M\in\hyperlink{LCset}{\mathcal{LC}}$, then $\omega_M$ satisfies \hyperlink{om5}{$(\omega_5)$} if and only if $\lim_{j\rightarrow+\infty}(M_j/j!)^{1/j}=+\infty$, see e.g. \cite[Lemma 2.8]{testfunctioncharacterization} and \cite[Lemma 2.4]{sectorialextensions} and the references mentioned in the proofs there.

Moreover, $\lim_{j\rightarrow+\infty}(M_j/j!)^{1/j}=+\infty$ if and only if $\lim_{j\rightarrow+\infty}\frac{\mu_j}{j}=+\infty$; see e.g. \cite[p. 104]{compositionpaper} and recall that by Stirling's formula $j\mapsto j!^{1/j}$ and $j\mapsto j$ grow comparable up to some positive constant.

For $M\in\hyperlink{LCset}{\mathcal{LC}}$ it is known that non-quasianalyticity implies $\lim_{j\rightarrow+\infty}\frac{\mu_j}{j}=+\infty$; see e.g. \cite[Prop. 4.4]{testfunctioncharacterization}. The growth of $j\mapsto\frac{\mu_j}{j}$ is crucial to see if the real-analytic functions are (strictly) contained in the corresponding ultradifferentiable class; more precisely $\lim_{j\rightarrow+\infty}\frac{\mu_j}{j}=+\infty$ if and only if the real-analytic functions are contained in the function classes under consideration, see e.g. again \cite[Prop. 2.12, p. 104]{compositionpaper}.\vspace{6pt}

If $0<s<1$, then for $\omega_s$ in Example \ref{promienentexample} corresponding to $\omega_{G^s}$ with $G^s:=(j!^s)_{j\in\NN}$ we get that $\lim_{j\rightarrow+\infty}(j!^s/j!)^{1/j}=0$. If $s=1$, then the limit is equal to $1$.\vspace{6pt}

Let now, more generally, the parameter $A>1$ be arbitrary (but fixed) and set $a_{j+1}:=2^{c_j}b_j$ for all $j\in\NN_{>0}$; i.e. $j$ in the exponent is replaced by a {\itshape fourth auxiliary sequence} of integers $(c_j)_{\in\NN_{>0}}$. Note that the conclusion of {\itshape Claim e} in {\itshape Step IV} in Theorem \ref{Abaninexamplethm} stays valid when replacing $j$ by any arbitrary non-decreasing $(c_j)_{j\in\NN_{>0}}$ satisfying $c_j\ge j$.

\begin{itemize}
\item[$(*)$] On each interval $I_{a_j,i}$ the mapping $k\mapsto\frac{\mu_k}{k}$ is clearly non-increasing.

\item[$(*)$] If $k=2^ia_j$ for some $j\in\NN_{>0}$ and $0\le i\le d_j-1$, i.e. positions where a jump occurs, then $\frac{\mu_k}{k}\le\frac{\mu_{k+1}}{k+1}\Leftrightarrow 1+\frac{1}{k}\le A$. This estimate is valid if $A\ge 2$ for all $k\in\NN_{>0}$ (as in Theorem \ref{Abaninexamplethm}) and, if $A>1$, for all $k$ sufficiently large depending on chosen $A$.

    Moreover, $\frac{\mu_{2^ia_j}}{2^ia_j}\le\frac{\mu_{2^{i+1}a_j}}{2^{i+1}a_j}\Leftrightarrow 2\le A$, for all $j\in\NN_{>0}$ and $0\le i\le d_j-1$.

\item[$(*)$] By \eqref{AbaninStepIVequ} and \eqref{AbaninStepIVequ1impl} we see
$$\frac{\mu_{b_j}}{b_j}=\left(\frac{A}{2}\right)^{d_j}\frac{\mu_{a_j}}{a_j},$$
and so $\frac{\mu_{b_j}}{b_j}<\frac{\mu_{a_j}}{a_j}$ if $A<2$ and $\frac{\mu_{b_j}}{b_j}\ge\frac{\mu_{a_j}}{a_j}$ if $A\ge 2$.

\item[$(*)$] In view of \eqref{AbaninStepIVequ4} we have for all $2^ib_j\le k<k+1\le 2^{i+1}b_j$, $j\in\NN_{>0}$ and $0\le i\le c_j-1$, that $\frac{\mu_k}{k}\ge\frac{\mu_{k+1}}{k+1}\Leftrightarrow 1+\frac{1}{k}\ge\sqrt{2}^{\frac{1}{2^ib_j}}\Leftrightarrow\left(1+\frac{1}{k}\right)^{2^{i+1}b_j}\ge 2$. Finally note that $k\le 2^{i+1}b_j$ and $k\ge b_1\ge 2a_1=2$ and so on each interval $I_{b_j,i}$ the mapping $k\mapsto\frac{\mu_k}{k}$ is non-increasing.

\item[$(*)$] If $k=2^{i+1}b_j$, $j\in\NN_{>0}$ and $0\le i\le c_j-2$, then $\frac{\mu_k}{k}\ge\frac{\mu_{k+1}}{k+1}\Leftrightarrow 1+\frac{1}{k}\ge\sqrt{2}^{\frac{1}{2^{i+1}b_j}}$ holds as shown before. If $k=2^{i+1}b_j$, $j\in\NN_{>0}$ and $i=c_j-1$ and so $k=a_{j+1}$, then $\frac{\mu_k}{k}\le\frac{\mu_{k+1}}{k+1}\Leftrightarrow 1+\frac{1}{k}\le A$ which holds for $A>1$ for all $k$ sufficiently large.

\item[$(*)$] By \eqref{AbaninStepIVequ4imp3} we get
$$\forall\;j\in\NN_{>0}:\;\;\;\frac{\mu_{a_{j+1}}}{a_{j+1}}=\frac{2^{\frac{c_j}{2}}}{2^{c_j}}\frac{\mu_{b_j}}{b_j}=\frac{1}{2^{\frac{c_j}{2}}}\frac{\mu_{b_j}}{b_j}<\frac{\mu_{b_j}}{b_j}.$$
\end{itemize}

Now we formulate and prove the second main statement concerning weights $\omega_M$.

\begin{theorem}\label{Abaninexamplethm1}
There exist (normalized) quasianalytic BMT-weight functions $\omega=\omega_M$ which are not equivalent to an almost sub-additive weight function.

In addition, we can achieve different situations concerning the growth of $k\mapsto\frac{\mu_k}{k}$:

\begin{itemize}
\item[$(*)$] Case $(a)$: $\lim_{k\rightarrow+\infty}\frac{\mu_k}{k}=0$.

\item[$(*)$] Case $(b)$: $\lim_{k\rightarrow+\infty}\frac{\mu_k}{k}=+\infty$.

\item[$(*)$] Case $(c)$: $\liminf_{k\rightarrow+\infty}\frac{\mu_k}{k}=0$ and $\limsup_{k\rightarrow+\infty}\frac{\mu_k}{k}=+\infty$.
\end{itemize}
\end{theorem}

\demo{Proof}
For all cases we use the analogous constructions and definitions as in Theorem \ref{Abaninexamplethm} but apply some technical modifications/generalizations and involve the comments on monotonicity before.\vspace{6pt}

{\itshape Case a:} We take the sequences $(a_j)_{j\in\NN_{>0}}$, $(b_j)_{j\in\NN_{>0}}$ and $(d_j)_{j\in\NN_{>0}}$ as considered in \eqref{AbaninStepIVequ} but then let us fix some $1<A_0<2$ and write $q:=\frac{2}{A_0}$. For $(d_j)_{j\in\NN_{>0}}$ we assume that the sequence is (strictly) increasing with $d_1\ge 1$. We also set $\mu_{a_1}=\mu_1:=1$ as ``initial condition'' instead of $\mu_{a_1}=2$ (which is required to prove \eqref{AbaninStepIVequ6} for $j=1$ but in the recent setting this property fails in any case as we are going to see). Next let us show that
\begin{equation}\label{Abaninexamplethmremequ}
\forall\;j\in\NN_{>0}\;\forall\;1<A\le A_0:\;\;\;\frac{a_j}{\mu_{a_j}}\frac{(2/A)^{d_j}-1}{(2/A)-1}\ge 1.
\end{equation}
First, since $\sum_{i=0}^kt^i=\frac{t^{k+1}-1}{t-1}$ for any $t>1$ we see that $t\mapsto\frac{t^{k+1}-1}{t-1}$ is non-decreasing on $(1,+\infty)$ for any fixed $k\in\NN$ and hence in \eqref{Abaninexamplethmremequ} it suffices to consider $q=\frac{2}{A_0}$. Then we use induction on $j$: The case $j=1$ gives $\frac{a_1}{\mu_{a_1}}\frac{q^{d_1}-1}{q-1}=\frac{q^{d_1}-1}{q-1}\ge 1$ which is clear. For the induction step we have
$$\frac{a_{j+1}}{\mu_{a_{j+1}}}=\frac{2^jb_j}{2^{\frac{j}{2}}\mu_{b_j}}=\frac{2^j2^{d_j}a_j}{2^{\frac{j}{2}}A^{d_j}\mu_{a_j}}=2^{\frac{j}{2}}\left(\frac{2}{A}\right)^{d_j}\frac{a_j}{\mu_{a_j}}\ge 2^{\frac{j}{2}}q^{d_j}\frac{a_j}{\mu_{a_j}}\ge 2^{\frac{j}{2}}q^{d_j}\frac{q-1}{q^{d_j}-1}\ge\frac{q-1}{q^{d_{j+1}}-1},$$
where the last desired estimate is equivalent to having $2^{\frac{j}{2}}(q^{d_{j+1}}-1)\ge 1-\frac{1}{q^{d_j}}$ and this clearly holds.
(However, note that for parameters $A\le 2$ the choice \eqref{AbaninStepIVequ5} and hence the estimate \eqref{AbaninStepIVequ6} is impossible.)

Then fix a real parameter $A$ with $1<A\le A_0$ and follow the proof of {\itshape Step IV.} All parts except {\itshape Claim b} follow by a word-by-word repetition of the given arguments there. Concerning {\itshape Claim a} note that normalization is still valid by $\mu_1=\mu_{a_1}=1$ and in {\itshape Claim d} we have $\liminf_{j\rightarrow+\infty}\frac{\mu_{4j}}{\mu_j}\ge\min\{2^{1/4},A\}>1$; so the value $2^{1/4}$ has to be replaced by $A$ if $A<2^{1/4}$ which is excluded by the choice $A>2$ in the proof before.

However, \hyperlink{nq}{$(\on{nq})$} fails: By \eqref{Abaninexamplethmremequ} we estimate as follows for all $j\in\NN_{>0}$:
\begin{align*}
\sum_{k=a_j+1}^{b_j}\frac{1}{\mu_k}&=\sum_{i=0}^{d_j-1}\sum_{k\in I_{a_j,i}}\frac{1}{\mu_k}=\sum_{i=0}^{d_j-1}\frac{2^ia_j}{A^{i+1}\mu_{a_j}}=\frac{a_j}{A\mu_{a_j}}\sum_{i=0}^{d_j-1}\left(\frac{2}{A}\right)^i=\frac{a_j}{A\mu_{a_j}}\frac{(2/A)^{d_j}-1}{(2/A)-1}\ge\frac{1}{A}.
\end{align*}
Since all appearing summands in the series under consideration are non-negative the previous estimate immediately verifies $\sum_{k\ge 1}\frac{1}{\mu_k}=+\infty$ as desired.

{\itshape Step V} follows again by varying the parameter $A\in(1,A_0]$; note that the choice of the numbers $d_j$ is not depending on given $A$ and so a direct comparison of the quotients is again possible. Finally, we get
$$\forall\;j\in\NN_{>0}:\;\;\;\frac{\mu_{a_{j+1}}}{a_{j+1}}=\frac{1}{2^{\frac{j}{2}}}\frac{\mu_{b_j}}{b_j}=\left(\frac{A}{2}\right)^{d_j}\frac{1}{2^{\frac{j}{2}}}\frac{\mu_{a_j}}{a_j},$$
and so by taking into account the comments just before the statement in this section it follows that $\lim_{k\rightarrow+\infty}\frac{\mu_k}{k}=0$ (for any parameter $A\in(1,A_0]$).\vspace{6pt}

{\itshape Case b:} First let us fix a value $A>2$. Second, we introduce {\itshape four} auxiliary sequences (of positive integers) $(a_j)_{j\in\NN_{>0}}$, $(b_j)_{j\in\NN_{>0}}$, $(c_j)_{j\in\NN_{>0}}$ and $(d_j)_{j\in\NN_{>0}}$ with $a_j<b_j<a_{j+1}$ for all $j\in\NN_{>0}$ and being defined as follows:
\begin{equation}\label{AbaninStepIVeququasi}
a_1:=1,\hspace{20pt}b_j:=2^{d_j}a_j,\hspace{20pt}a_{j+1}:=2^{c_j}b_j,\hspace{20pt}j\in\NN_{>0}.
\end{equation}
Thus $2^j$ in \eqref{AbaninStepIVequ} is replaced by $2^{c_j}$ and so $(b_j,a_{j+1}]$, $j\in\NN_{>0}$, is now split into $c_j$-many subintervals $I_{b_j,i}:=(2^ib_j,2^{i+1}b_j]$, $j\in\NN_{>0}$, $i\in\NN$ with $0\le i\le c_j-1$.

The definition of $\mu$ is then analogous according to \eqref{AbaninStepIVequ1} and \eqref{AbaninStepIVequ4} and again we consider the ``initial condition'' $\mu_{a_1}=\mu_1:=1$. Hence
\begin{equation}\label{AbaninStepIVequ1implquasi}
\forall\;j\in\NN_{>0}:\;\;\;\mu_{b_j}=A^{d_j}\mu_{a_j},\hspace{30pt}\mu_{a_{j+1}}=2^{\frac{c_j}{2}}\mu_{b_j},
\end{equation}
meaning that \eqref{AbaninStepIVequ1impl} remains unchanged and in \eqref{AbaninStepIVequ4imp3} the power $j$ is replaced by $c_j$. However, the crucial difference compared with the construction in Theorem \ref{Abaninexamplethm} is now that we choose $(d_j)_{j\in\NN_{>0}}$, $(c_j)_{j\in\NN_{>0}}$ iteratively as follows: Both sequences are required to be strictly increasing with
\begin{equation}\label{AbaninStepIVequ5quasibase}
d_1\ge 1,\hspace{30pt}c_j\ge j,\;\;\;\forall\;j\in\NN_{>0},
\end{equation}
such that
\begin{equation}\label{AbaninStepIVequ5quasi}
\forall\;j\in\NN_{>0}:\;\;\;\frac{\mu_{b_j}}{b_j}\ge\sqrt{2}^{c_j-1}\log(j+1),
\end{equation}
and
\begin{equation}\label{AbaninStepIVequ5quasi1}
\forall\;j\in\NN_{>0}:\;\;\;\sqrt{2}^{c_j-1}\ge\frac{\mu_{b_j}}{b_j}\frac{1}{j\log(j+1)}+1.
\end{equation}
These choices are possible: Recall that $A>2$ and $\frac{\mu_{b_j}}{b_j}=\frac{\mu_{a_j}}{a_j}\left(\frac{A}{2}\right)^{d_j}$ for all $j\in\NN_{>0}$. Hence, when given $\mu_{a_j}$ and $a_j$, then the difference $\frac{\mu_{b_j}}{b_j}\frac{1}{\log(j+1)}-\frac{\mu_{b_j}}{b_j}\frac{1}{j\log(j+1)}-1$ becomes as large as desired when $d_j$ is chosen sufficiently large. Since $c_j$ is only required to define $\mu_{a_{j+1}}$, $a_{j+1}$, and so not depending on $\mu_{b_j}$, $b_j$, both \eqref{AbaninStepIVequ5quasi} and \eqref{AbaninStepIVequ5quasi1} can be guaranteed.\vspace{6pt}

All properties from {\itshape Step IV} except {\itshape Claim b} dealing with \hyperlink{nq}{$(\on{nq})$} follow analogously as shown in the claims there; for {\itshape Claim e} note that $c_j\ge j$ for all $j$. Concerning the failure of \hyperlink{nq}{$(\on{nq})$} by definition, \eqref{AbaninStepIVequ5quasi1} (and recall \eqref{AbaninStepIVequ4imp2}) we estimate as follows for all $j\in\NN_{>0}$:
\begin{align*}
\sum_{k=b_j+1}^{a_{j+1}}\frac{1}{\mu_k}&=\sum_{i=0}^{c_j-1}\sum_{k\in I_{b_j,i}}\frac{1}{\mu_k}\ge\sum_{i=0}^{c_j-1}\frac{2^{i+1}b_j-2^ib_j}{\mu_{2^{i+1}b_j}}=\sum_{i=0}^{c_j-1}\frac{2^ib_j}{2^{\frac{i+1}{2}}\mu_{b_j}}
\\&
=\frac{1}{\sqrt{2}}\frac{b_j}{\mu_{b_j}}\sum_{i=0}^{c_j-1}2^{\frac{i}{2}}=\frac{1}{\sqrt{2}}\frac{b_j}{\mu_{b_j}}\frac{\sqrt{2}^{c_j}-1}{\sqrt{2}-1}\ge\frac{1}{\sqrt{2}}\frac{b_j}{\mu_{b_j}}\frac{\sqrt{2}^{c_j}-\sqrt{2}}{\sqrt{2}-1}
\\&
=\frac{b_j}{\mu_{b_j}}\frac{\sqrt{2}^{c_j-1}-1}{\sqrt{2}-1}\ge\frac{1}{\sqrt{2}-1}\frac{1}{j\log(j+1)}.
\end{align*}
Since all appearing summands in the series under consideration are non-negative the previous estimate immediately verifies $\sum_{k\ge 1}\frac{1}{\mu_k}=+\infty$ as desired.

Let us see that here we get $\lim_{k\rightarrow+\infty}\frac{\mu_k}{k}=+\infty$ as well: By definition $\frac{\mu_{b_j}}{b_j}=\frac{\mu_{a_j}}{a_j}\left(\frac{A}{2}\right)^{d_j}\rightarrow+\infty$ as $j\rightarrow+\infty$, whereas \eqref{AbaninStepIVeququasi}, \eqref{AbaninStepIVequ1implquasi} and \eqref{AbaninStepIVequ5quasi} yield $$\frac{\mu_{a_{j+1}}}{a_{j+1}}=\frac{2^{\frac{c_j}{2}}\mu_{b_j}}{2^{c_j}b_j}=\frac{1}{2^{\frac{c_j}{2}}}\frac{\mu_{b_j}}{b_j}\ge\frac{1}{\sqrt{2}}\log(j+1)\rightarrow+\infty$$
as $j\rightarrow+\infty$. By the comments made at the beginning of this section we are done.\vspace{6pt}

{\itshape Case c:} Here we take the same construction and definition as in {\itshape Case b} but for the sequences $(d_j)_{j\in\NN_{>0}}$ and $(c_j)_{j\in\NN_{>0}}$, we assume the following: Again we require \eqref{AbaninStepIVequ5quasibase} but instead of \eqref{AbaninStepIVequ5quasi} and \eqref{AbaninStepIVequ5quasi1} iteratively we choose $(d_j)_{j\in\NN_{>0}}$ and $(c_j)_{j\in\NN_{>0}}$ to grow fast enough to ensure
\begin{equation}\label{sequenceccasec0}
\forall\;j\in\NN_{>0}:\;\;\;\frac{\mu_{b_j}}{b_j}=\frac{\mu_{a_j}}{a_j}\left(\frac{A}{2}\right)^{d_j}>j,
\end{equation}
and
\begin{equation}\label{sequenceccasec}
\forall\;j\in\NN_{>0}:\;\;\;\sqrt{2}^{c_j}\ge j\frac{\mu_{b_j}}{b_j}.
\end{equation}
\eqref{sequenceccasec0} implies $\limsup_{k\rightarrow+\infty}\frac{\mu_k}{k}=+\infty$ and \eqref{sequenceccasec} yields, by taking into account \eqref{AbaninStepIVeququasi} and \eqref{AbaninStepIVequ1implquasi}, that $$\frac{\mu_{a_{j+1}}}{a_{j+1}}=\frac{1}{2^{\frac{c_j}{2}}}\frac{\mu_{b_j}}{b_j}\le\frac{1}{j},$$
and so  $\liminf_{k\rightarrow+\infty}\frac{\mu_k}{k}=0$ follows.

Now let us see that \eqref{sequenceccasec} already implies \eqref{AbaninStepIVequ5quasi1} for all $j$ large and this is sufficient to verify the quasianalyticity as in {\itshape Case b:} We have $\sqrt{2}^{c_j-1}-1\ge\frac{j}{\sqrt{2}}\frac{\mu_{b_j}}{b_j}-1$ and clearly, in view of \eqref{sequenceccasec0},  $\frac{j}{\sqrt{2}}\frac{\mu_{b_j}}{b_j}-1\ge\frac{\mu_{b_j}}{b_j}\frac{1}{j\log(j+1)}\Leftrightarrow\frac{\mu_{b_j}}{b_j}\left(\frac{j}{\sqrt{2}}-\frac{1}{j\log(j+1)}\right)\ge 1$ for all $j$ large enough.
\qed\enddemo

\bibliographystyle{plain}
\bibliography{Bibliography}
\end{document}